\documentclass[11pt]{article}

\usepackage{amsfonts,amssymb}
\usepackage{diagrams}

\usepackage[latin1]{inputenc}  
\usepackage[T1]{fontenc}

\def\pplus{\tilde{+}}

\newcounter{paragrafsubsub}[subsubsection]
\renewcommand{\theparagrafsubsub}{%
\thesubsubsection.\roman{paragrafsubsub}}
\newcommand{\paragrafsubsub}{%
\refstepcounter{paragrafsubsub}
{\bf \theparagrafsubsub}\hspace{0.2em}--- }

\newcounter{paragrafsub}[subsection]
\renewcommand{\theparagrafsub}{\thesubsection.\arabic{paragrafsub}}
\newcommand{\paragrafsub}{%
\refstepcounter{paragrafsub}
{\bf \theparagrafsub}\hspace{0.2em}--- }

\newcounter{paragraf}[section]
\renewcommand{\theparagraf}{\thesection.\arabic{paragraf}}
\newcommand{\paragraf}{%
\refstepcounter{paragraf}
{\bf \theparagraf}\hspace{0.2em}--- }

\newcommand\paragraphe{%
\par \indent
\ifcase\value{subsection} %
\paragraf
\else
\ifcase\value{subsubsection}\paragrafsub %
\else\paragrafsubsub
\fi\fi
}

\def\longto{\longrightarrow}
\def\quot{/\hspace{-.5ex}/}

\def\Id{{\rm Id}}

\def\Rep{{\rm Rep}}
\def\Pic{{\rm Pic}}
\def\Hom{{\rm Hom}}
\def\GL{{\rm GL}}\def\SL{{\rm SL}}

 \def\Hcal{{\mathcal H}}\def\Wcal{{\mathcal W}}

\def\Orb{{\mathcal O}}
\def\PP{{\mathbb P}}\def\NN{{\mathbb N}}
\def\QQ{{\mathbb Q}}\def\ZZ{{\mathbb Z}}
\def\CC{{\mathbb C}}

\def\kk{{\mathbb K}}

\def\GL{{\rm GL}}

\def\tc{{\mathcal{TC}}}
\def\ac{{\mathcal{AC}}}
\def\Pol{{\mathcal{P}}}
\def\Face{{\mathcal F}}
\def\Li{{\mathcal{L}}}\def\Mi{{\mathcal{M}}}

\newtheorem{lemma}{Lemma}
\newtheorem{prop}{Proposition}
\newtheorem{theo}{Theorem}

\newenvironment{proof}{{\noindent\bf Proof.}}{\hfill $\square$}
\newenvironment{defin}{{\noindent\bf Definition.}}{\\}
\newenvironment{remark}{{\noindent\bf Remark.}}{}

\begin{document}
\title{GIT-cones and quivers}
\author{N. Ressayre\footnote{Universit{\'e} Montpellier II,
D{\'e}partement de Math{\'e}matiques,
Case courrier 051-Place Eug{\`e}ne Bataillon,
34095 Montpellier Cedex 5, France. {\tt ressayre@math.univ-montp2.fr}}}

\maketitle
\begin{abstract}
In this work, we improve results of \cite{GITEigen,GITEigen2} about GIT-cones associated to the
action of any reductive group $G$ on a projective variety $X$.
These results are applied to give a short proof of a Derksen-Weyman's Theorem which parametrizes
bijectively the faces of a rational cone associated to any quiver without oriented cycle. 
\end{abstract}

\tableofcontents 

\newpage
\section{Introduction}

We work over the complex numbers field $\CC$.
Let $Q=(Q_0,Q_1)$ be a quiver without oriented cycle. 
Here, $Q_0$ is the set of vertices and $Q_1$ the set of arrows.
Let $\beta=(\beta(s))_{s\in Q_0}$  be a vector dimension of $Q$ and
$\Rep(Q,\beta)$ the vector space of the representations of dimension vector $\beta$.
The group $\GL(\beta)=\prod_{s\in Q_0}\GL(\beta(s))$  acts naturally on $\Rep(Q,\beta)$.
We consider the group $\Gamma$ of the characters of $\GL(\beta)$; it is isomorphic to $\ZZ^{Q_0}$.
We consider the cone $\Sigma(Q,\beta)$ in $\Gamma\otimes\QQ$ generated the elements  $\sigma\in\Gamma$ 
such that there exists a non zero regular function $f\in\CC[\Rep(Q,\beta)]$ such that $g.f=\sigma(g)f$
for any $g\in\GL(\beta)$.
Actually, it is a convex polyhedral cone.
In \cite{DW:saturation,DW:comb}, Derksen-Weyman showed that Horn's cones can be obtained in such a way.
This is an important motivation for the study of these cones.
Here, we use general methods of Geometric Invariant Theory to give a proof of a Derksen-Weyman's Theorem (see 
\cite{DW:comb}) which parametrizes bijectively the faces of $\Sigma(Q,\beta)$.\\
 
In Section~\ref{sec:gen}, we improve results of \cite{GITEigen} about GIT-cones in general. 
In particular, Theorem~\ref{th:wc2face} is an improvement of \cite[Theorem~7]{GITEigen},  
and Theorem~\ref{th:faceL} of \cite[Theorem~6]{GITEigen}.
So, the examples  $\Sigma(Q,\beta)$ enlight the general theory of GIT-cones.\\

The Horn's cones can also be obtained as GIT-cones for the action of the linear groups on products of
complete flag varieties. This point of view was used \cite{sjamaar,BK,GITEigen,GITEigen2} for example.
Whereas, in the literature this GIT approach of Horn's problem was distinct from the quiver one,
this work prove that one can use essentially the same techniques in the two cases.

\section{Well covering pairs and GIT-cones}
\label{sec:gen}

\subsection{Well covering pairs}
\label{sec:def}

Let $G$ be a reductive group acting on a smooth variety $X$.
Let $\lambda$ be a one parameter subgroup of $G$.
Let $G^\lambda$ denote the centralizer of $\lambda$ in $G$.
We consider the usual 
parabolic subgroup $P(\lambda)$ associated to $\lambda$ with Levi subgroup $G^\lambda$:
$$
P(\lambda)=\left \{
g\in G \::\:
\lim_{t\to 0}\lambda(t).g.\lambda(t)^{-1} 
\mbox{  exists in } G \right \}.
$$

Let $C$ be an irreducible  component of the fix point set $X^\lambda$ of $\lambda$ in $X$.
We also consider the Bialinicky-Birula cell $C^+$ associated to $C$:
$$
C^+=\{x\in X\ |\ \lim_{t\to 0}\lambda(t)x\in C\}.
$$
Then, $C$ is stable by the action of $G^\lambda$ and $C^+$ by the action of $P(\lambda)$.

Consider over $G\times C^+$ the action of $G\times P(\lambda)$ given by the formula 
(with obvious notation): $(g,p).(g',y)=(gg'p^{-1},py)$.
Consider the  quotient $G\times_{P(\lambda)}C^+$ of $G\times C^+$ by the action of 
$\{e\}\times P(\lambda)$.
The class of a pair $(g,y)\in G\times C^+$ in $G\times_{P(\lambda)} C^+$ is denoted by $[g:y]$.

The action of $G\times\{e\}$ induces an action of $G$ on $G\times_{P(\lambda)}C^+$.
Moreover, the first projection $G\times C^+\longto G$ induces a $G$-equivariant map  
$\pi\,:\,G\times_{P(\lambda)} C^+\longto G/P(\lambda)$ which is a locally trivial fibration 
with fiber $C^+$. 
Consider also the $G$-equivariant map 
$$\eta\,:\,G\times_{P(\lambda)}C^+\longto X,\,[g:y]\mapsto gy.$$

\begin{defin}
The pair $(C,\lambda)$ is said to be {\it dominant} if $\eta$ is.
 The pair $(C,\lambda)$ is said to be {\it well covering} if $\eta$ induces an isomorphism over
an open subset of $X$ intersecting $C$. 
\end{defin}

Let $\Li\in\Pic^G(X)$.
Let $x$ be any point in $C$.
Since $\lambda$ fixes $x$, it induces a linear action of 
the group $\kk^*$ on the fiber $\Li_x$.
This action defines a character of $\kk^*$, that is, an element of
$\ZZ$ denoted by $\mu^\Li(x,\lambda)$. 
One easily checks that $\mu^\Li(x,\lambda)$ does not depends on $x\in C$; it will be denoted by
$\mu^\Li(C,\lambda)$. 

\subsection{Total cones and well covering pair}

\paragraphe
Consider the convex cones $\tc^G(X)$ generated in $\Pic^G(X)_\QQ$ by the $\Li$'s  in $\Pic^G(X)$ 
which have non zero $G$-invariant sections. 
We will denote by $X^{\rm ss}(\Li)$ the open subset of the $x$'s in $X$ such that 
some positive integer $n$, there exists a $G$-invariant section of $\Li^{\otimes n}$ 
such that $\sigma(x)\neq 0$. Note that this definition is standard, only if $\Li$ 
is ample. 
Since $X^{\rm ss}(\Li)=X^{\rm ss}(\Li^{\otimes n})$ 
(for any positive integer $n$), one can define $X^{\rm ss}(\Li)$ if $\Li\in\Pic^G(X)_\QQ$.

Let $(C,\lambda)$ be a dominant  pair. Since $\Li\mapsto \mu^\Li(C,\lambda)$ is a group morphism,
it induces a linear map from $\Pic^G(X)_\QQ$ to $\QQ$, also denoted by $\mu^\Li(C,\lambda)$.
By \cite[Lemma~7]{GITEigen}, $\tc^G(X)$ is contained in the halfspace $\mu^\Li(C,\lambda)\leq 0$.
In particular, intersecting $\tc^G(X)$ with the hyperplane $\mu^\Li(C,\lambda)=0$, one obtains 
a face $\Face(C)$ of $\tc^G(X)$.
Indeed, the following lemma shows that the face only depends on $C$:

\begin{lemma}
\label{lem:FC}
Let $(C,\lambda)$ be a dominant pair.
Then, $\Face(C)$ is the set of $\Li\in\Pic^{G}(X)_\QQ$ such that $X^{\rm ss}(\Li)$ 
intersects $C$. 
\end{lemma}

\begin{proof}
Assume that   $X^{\rm ss}(\Li)$  intersects $C$. 
Let $\sigma$ be a $G$-invariant section of $\Li^{\otimes n}$ and $x\in C$ such that 
$\sigma(x)\neq 0$.
We have: $\lambda(t)\sigma(x)$ is equal on one hand to $t^{n\mu^{\Li}(C,\lambda)}\sigma(x)$ and on the hand to
$\sigma(\lambda(t)x)=\sigma(x)$. It follows that $ \mu^{\Li}(C,\lambda)=0$.

Conversely, let $\Li\in\Face(C)$. Let $\sigma$ be a non zero $G$-invariant section of $\Li^{\otimes n}$.
Since $\eta$ is dominant, $\sigma$ is non identically zero on $C^+$. 
By \cite[Proposition~5]{GITEigen}, $\mu^{n\Li}(C,\lambda)=0$ implies that
 $\sigma$ is non identically zero on $C$.
In particular, $X^{\rm ss}(\Li)$ intersects $C$. 
\end{proof}\\

\begin{remark}
  In the above proof, we have only used that $\eta$ is dominant.
\end{remark}

The following theorem is an improvement of \cite[Theorem~7]{GITEigen}.

\begin{theo}
\label{th:wc2face}
We assume the rank of  $\Pic^G(X)$ is finite and consider $\tc^G(X)$.
Let $(C,\lambda)$ be a well covering pair.

The rank of $\Pic^{G^\lambda}(C)$ is finite.
The codimension of $\Face(C)$ in $\Pic^G(X)_\QQ$ equals the codimension of 
$\tc^{G^\lambda}(C)$ in $\Pic^{G^\lambda}(C)_\QQ$.
More precisely, 
the restriction morphism induces an isomorphism
from $\Pic^G(X)_\QQ/\langle \Face(C)\rangle$ onto
 $\Pic^{G^\lambda}(C)_\QQ/\langle \tc^{G^\lambda}(C)\rangle$.
\end{theo}

\begin{proof}
  Let $\Omega$ be a $G$-stable open subset of $X$ such that the natural map 
$G\times_{P(\lambda)}(C^+\cap \Omega)\longto \Omega$ is an isomorphism.
Since $(C,\lambda)$ is well covering one can find such an $\Omega$ intersecting $C$.
By \cite[Lemma~1]{GITEigen}, $\Pic^G(\Omega)$ is isomorphic to $\Pic^{G^\lambda}(C\cap \Omega)$.

Let $E_1,\cdots,E_s$ (resp. $D_1,\cdots,D_t$) be the irreducible components of codimension one of 
$X-\Omega$ (resp. $C-\Omega$). 
Since $G$ and $G^\lambda$ are connected the $E_i$'s and the $D_i$'s are respectively $G$ and 
$G^\lambda$-stable. We consider the associated $G$ and $G^\lambda$-linearized line bundles $\Li_{E_i}$
and $\Li_{D_i}$. 

Consider the following diagram:

\begin{diagram}
\oplus_i\QQ\Li_{E_i}&\rTo&\Pic^G(X)_\QQ \\
                    &    &             &\rdTo^{\pi_X}\\
                    &    & \dTo        &\ \ \ \ \ \  &\Pic^G(\Omega)_\QQ\simeq\Pic^{G^\lambda}(C\cap \Omega)_\QQ\\
                    &    &             &\ruTo_{\pi_C}\\
\oplus_i\QQ\Li_{D_i}&\rTo&\Pic^{G^\lambda}(C)_\QQ
\end{diagram}

Since $X$ and so $C$ are  smooth, the maps $\pi_X$ and $\pi_C$ are surjective. 
By construction, $\Li_{D_i}$ belongs to $\tc^{G^\lambda}(C)$.
Moreover, each $\Li_{E_i}$ has a $G$-invariant section with $E_i$ has zero locus. 
Since $E_i$ does not contains $C$, this proves that $\Li_{E_i}\in\Face(C)$.
So, it is sufficient to prove that $\pi_X(\Face(C))=\pi_C(\tc^{G^\lambda}(C))$. 

Let $\Li\in\Face(C)$.
Since $X^{\rm ss}(\Li)$ intersects $C$, $\Li_{|C}$ belongs to $\tc^{G^\lambda}(C)$.
So, $\pi_X(\Face(C))\subset \pi_C(\tc^{G^\lambda}(C))$. 
 
Conversely, let $\Li$ be a $G^\lambda$-linearized line bundle on $C$ which belongs to 
$\tc^{G^\lambda}(C)$ and $\sigma$ be a non zero $G^\lambda$-invariant section of 
$\Li$. Let $\tilde\Li$ be the $G$-linearized line bundle on $\Omega$ associated to $\pi_C(\Li)$, and
$\tilde\sigma$ the $G$-invariant section of $\tilde\Li$ associated to $\sigma$.
Now, let $\Mi\in\Pic^G(X)$ such that $\pi_X(\Mi)=\tilde\Li$.
The section $\tilde\sigma$ induces a non zero $G$-invariant rational section of $\Mi$, and so a
non zero regular $G$-invariant section $\sigma'$ of $\Mi'=\Mi+\sum_i \Li_{E_i}^{\otimes n_i}$ for some non negative
integers $n_i$. Since no $E_i$ contains $C$, $\sigma'$ is not identically zero on $C$; in particular,
$\Mi'\in \Face(C)$.
Since, $\pi_X(\Mi')=\pi_X(\Mi)=\pi_C(\Li)$, it follows that $\pi_X(\Face(C))\supset \pi_C(\tc^{G^\lambda}(C))$.  
Note that details about the above argue can be found in the proof of \cite[Theorem~7]{GITEigen}.
\end{proof}\\

\paragraphe
The set of ample $G$-linearized line bundles generate an open convex $\Pic^G(X)^+_\QQ$ in $\Pic^G(X)_\QQ$.
We set: $\ac^G(X)=\Pic^G(X)^+_\QQ\cap\tc^G(X)$.
If $\Face$ is a face of $\tc^G(X)$, $\Face^0$ denotes its intersection with $\Pic^G(X)^+_\QQ$.
To any ample $\Li\in\Pic^G(X)_\QQ$ which does not belong to $\ac^G(X)$, using mainly 
Kempf's Theorem, we associated in \cite{GITEigen} a well covering pair $(C,\lambda)$
(actually, defined up to conjugacy).
Note that, to do this we need to fix a ``norm'' on the set of one parameter subgroups of $G$ invariant 
by conjugacy.  When restricted to $Y(S)\otimes \QQ$, (where $Y(S)$ denotes the group of one parameter 
subgroups  of any subtorus $S$ of $G$), this norm becomes the norm associated to a scalar product.
The face $\Face(C)$ is also denoted by $\Face(\Li)$. 
\cite[Theorem~6]{GITEigen} asserts that any face of $\ac^G(X)$ equals $\Face^\circ(\Li)$ for some 
ample $\Li\not\in\ac^G(X)$. 
Here, we need an improvement of this result.

Let $K^\lambda$ denote the neutral component of the kernel of the action of $G^\lambda$ on $C$.
Note that $\Pic^{G^\lambda/K^\lambda}(C)_\QQ$ is naturally embedded in $\Pic^{G^\lambda}(C)_\QQ$;
and contains $\tc^{G^\lambda}(C)$.

\begin{theo}
\label{th:faceL}
Let $\Face$ be a face of $\ac^G(X)$.
Consider the set $\Delta(\Face)$ of ample $\Li\not\in\ac^G(X)$ such that $\Face=\Face^\circ(\Li)$.

Then,    
\begin{enumerate}
\item There exists $\Li\in \Delta(\Face)$ such that the associated pair $(C,\lambda)$ satisfies
$\ac^{G^\lambda}(C)$ has non empty interior in $\Pic^{G^\lambda/K^\lambda}(C)$.
\item $\Delta(\Face)$ contains a non empty open subset of $\Pic^G(X)_\QQ^+-\ac^G(X)$.
\end{enumerate}

\end{theo}

\begin{proof}
  Let $S$ be a maximal torus of $K^\lambda$.
Let $S'$ be a torus of $G^\lambda$ such that the product induces an isogeny from $S\times S'$ onto 
a maximal torus $T$
of $G^\lambda$ and $Y(S)$ is orthogonal with $Y(S')$. 
Let $H^\lambda$ be the connected subgroup of $G^\lambda$ such that 
the product induces an isogeny $K^\lambda\times H^\lambda\longto G^\lambda$ 
and containing $S'$. 

Consider the restriction morphism $p\,:\,\Pic^G(X)\longto\Pic^{H^\lambda}(C)$.
Note also that $\Pic^{H^\lambda}(C)_\QQ$ is isomorphic to $\Pic^{G^\lambda/K^\lambda}(C)_\QQ$.

Let $\Li\in\Delta(\Face)$.
We assume that the interior of $\ac^{G^\lambda}(C)$ in $\Pic^{G^\lambda/K^\lambda}(C)$ is empty.
By \cite[Lemma~9]{GITEigen}, $p(\Li)$ belongs to $\ac^{H^\lambda}(C)$.
By Theorem~\ref{th:wc2face}, one can find a neighbor $\Li_\epsilon$ of $\Li$ such that $p(\Li_\epsilon)$
does not belong to the span of $\ac^{G^\lambda}(C)$. By \cite[Lemma~10]{GITEigen}, on can find such a
$\Li_\epsilon$ such that the face of $\ac^{H^\lambda}(C)$ viewed from $\Li_\epsilon$ is the whole 
$\ac^{H^\lambda}(C)$. By \cite[Lemma~11]{GITEigen}, one may assume that $\Face(\Li_\epsilon)=\Face$.
Moreover the proof of \cite[Lemma~11]{GITEigen} shows that $C_\epsilon$ (with obvious notation) is
strictly contained in $C$.
By induction on the dimension of $C$ we just proved that there exists $\Li\in\Delta(\Face)$ such that
 interior of $\ac^{G^\lambda}(C)$ in $\Pic^{G^\lambda/K^\lambda}(C)$ is non empty.\\

By the same argue as above, one can prove that there exists $\Li\in\Delta(\Face)$
such that $p(\Li)$ belongs to the interior of
$\ac^{H^\lambda}(C)\cap {\rm Im}p$ in ${\rm Im}p$.
In this case, $p^{-1}(\ac^{H^\lambda}(C))\cap (\Pic^G(X)^+_\QQ-\ac^G(X))$ has non empty interior in 
$\Pic^G(X)_\QQ$. Moreover, by \cite[Lemma~10]{GITEigen}, this set is contained in $\Delta(\Face)$.
The second assertion is proved.
\end{proof}

\section{Application to quiver representations}

\subsection{Definitions}
\label{sec:defcarquois}

In this section, we fix some classical notation about quiver representations.

Let $Q$ be a quiver (that is, a finite oriented graph) with vertexes $Q_0$ and arrows $Q_1$.
An arrow $a\in Q_1$ has initial vertex $ia$ and terminal one $ta$.
A representation $R$ of $Q$ is a family $(V(s))_{s\in Q_0}$ of finite dimensional vector spaces and 
a family of linear maps $u(a)\in {\rm Hom}(V(ia),V(ta))$ indexed by $a\in Q_1$.
The dimension vector of $R$ is the family $(\dim(V(s)))_{s\in Q_0}\in \NN^{Q_0}$.

Let us fix $\alpha\in \NN^{Q_0}$ and  a vector space $V(s)$ of dimension $\alpha(s)$ for each $s\in Q_0$. 
Set
$$
\Rep(Q,\alpha)=\bigoplus_{a\in Q_1}{\rm Hom}(V(ia),V(ta)).
$$
Consider also the groups:
$$
\GL(\alpha)=\prod_{s\in Q_0}\GL(V(s))
{\rm\ and\ }
\SL(\alpha)=\prod_{s\in Q_0}\SL(V(s).
$$
They acts naturally on $\Rep(Q,\alpha)$.

The character group of $\GL(\alpha)$ identifies with $\Gamma=\ZZ^{Q_0}$; to 
$\sigma\in \ZZ^{Q_0}$, we associate the character $\chi_\sigma$ defined by
 $\chi_\sigma\left( (g(s))_{s\in Q_0})\right)=\prod_{s\in Q_0}\det(g(s))^{\sigma(s)}$.

\subsection{Three cones}

\paragraphe
Consider the algebra $\CC[\Rep(Q,\alpha)]$ of regular functions on $\Rep(Q,\alpha)$ endowed with the $\GL(\alpha)$-action.
For $\sigma\in \ZZ^{Q_0}$, we denote $\CC[\Rep(Q,\alpha)]_\sigma$ the set of $f\in\CC[\Rep(Q,\alpha)]$ such that
for all $g\in\GL(\alpha)$, $g.f=\chi_\sigma(g)f$.
We embed $\Gamma=\ZZ^{Q_0}$ in $\Gamma_\QQ:=\QQ^{Q_0}$.
Let $\Sigma(Q,\alpha)$ denote the convex cone of $\Gamma_\QQ$ generated by the points 
$\sigma\in\Gamma$ such that $\CC[\Rep(Q,\alpha)]_{-\sigma}$ is non reduced to $\{0\}$.

\paragraphe
Consider the projective space $X=\PP(\Rep(Q,\alpha)\oplus \CC)$. The formula 
$$
g.(R,t)=(gR,t)\ \ \ \forall g\in\GL(\alpha)\,,R\in\Rep(Q,\alpha){\rm\ and\ }
t\in\CC,
$$
defines an action of $\GL(\alpha)$ on $X$ and a $\GL(\alpha)$-linearization $\Li_0\in\Pic^{\GL(\alpha)}(X)$ of the line
bundle ${\mathcal O}(1)$ on $X$.
We are now interested in the GIT-cone $\ac^{\GL(\alpha)}(X)$.
Since any line bundle on $X$ admitting non zero sections is ample, $\ac^{\GL(\alpha)}(X)=\tc^{\GL(\alpha)}(X)$.

For $n\in \ZZ$ and $\sigma\in \Gamma$, set $\Li(n,\sigma)=\Li_0\otimes\sigma\in\Pic^{\GL(\alpha)}(X)$.
Note that $\Li(n,\sigma)={\mathcal O}(n)$ as a   line bundle.
We have the following obvious

\begin{lemma}
\label{lem:PicGX}
  The map $\ZZ\times\Gamma\longto\Pic^{\GL(\alpha)}(X),\,(n,\sigma)\mapsto\Li(n,\sigma)$ is an isomorphism of groups.
Moreover, $\Li(n,\sigma)$ is ample if and only if $n$ is positive.
\end{lemma}

Lemma~\ref{lem:PicGX} allows to embed $\ac^{\GL(\alpha)}(X)$ in $\QQ\times\Gamma_\QQ$. 
Set $\Pol(Q,\alpha)=\ac^{\GL(\alpha)}(X)\cap \{1\}\times\Gamma_\QQ$.
General properties of $\ac^{\GL(\alpha)}(X)$ imply that $\Pol(Q,\alpha)$ is  closed convex rational and polyhedral.
Moreover, it is contained in the convex hull of the restrictions to the center of $\GL(\alpha)$ of the weights
of a maximal torus of $\GL(\alpha)$ on $\Rep(Q,\alpha)\oplus\CC$; and so, it is compact. 
Finally,  $\Pol(Q,\alpha)$ is a rational polytope in $\Gamma_\QQ$.
Moreover, $\ac^{\GL(\alpha)}(X)$ is the convex cone generated by $\Pol(Q,\alpha)$ in such a  way the faces of
$\ac^{\GL(\alpha)}(X)$ and $\Pol(Q,\alpha)$ correspond bijectively.
\\

\paragraphe
We consider $\Rep(Q,\alpha)$ as an open subset of $X$ by $R\mapsto (R,1)$; and we identify the complement
with $\PP(\Rep(Q,\alpha))$.

\begin{prop}
\label{prop:0sommet}
\begin{enumerate}
\item We have: $X^{\rm ss}(\Li_0)=\Rep(Q,\alpha)$.
\item The point $0\in\Gamma_\QQ$ is a vertex of $\Pol(Q,\alpha)$.
\item The cone of $\Gamma_\QQ$ generated by $\Pol(Q,\alpha)$ is $\Sigma(Q,\alpha)$.
\end{enumerate}
\end{prop}

\begin{proof}
Since $Q$ has no oriented cycle, one can chose a numeration of the vertices such that 
the index of $ta$ is greater than the index of $ia$ for all $a\in Q_1$.
Consider the one parameter subgroup $\lambda_0$ of $\GL(\alpha)$ acting on the vector space 
corresponding to the vertex indexed by $i$ as an homothetie of coefficient $t^i$.

The point $0\in\Rep(Q,\alpha)\subset X$ is an isolated fixed point of $\lambda_0$.
Set $C_0=\{0\}$. One easily checks that $C_0^+=\Rep(Q,\alpha)$ and that $\lambda_0$ is central in
$\GL(\alpha)$: it follows that $(C_0,\lambda_0)$ is a well covering pair.
Let $\Face(C_0)$ (resp. $\Pol(C_0)$) denote the face of 
$\tc^{\GL(\alpha)}(X)$ (resp. $\Pol(Q,\alpha)$) associated to $(C_0,\lambda_0)$.

Viewed the action of the center of $\GL(\alpha)$ on the fiber over $0$,
 $\Face(C_0)$ is contained in $\QQ^+\Li_0$.

  Since $(R,t)\mapsto t$ is a $\GL(\alpha)$-invariant section of $\Li_0$, 
$\Rep(Q,\alpha)\subset X^{\rm ss}(\Li_0)$.
Then, $\Face(C_0)=\QQ^+\Li_0$.

Since $0$ is the only closed $\GL(\alpha)$-orbit in $\Rep(Q,\alpha)$,
$X^{\rm ss}(\Li_0)\quot \GL(\alpha)$ is only one point.
So, $X^{\rm ss}(\Li_0)$ contains only one closed orbit $\Orb$ which is contained in the closure 
$\GL(\alpha).0$. We deduce that $\Orb=\{0\}$ and that $X^{\rm ss}(\Li_0)=\Rep(Q,\alpha)$.

The last assertion of the proposition is a direct application of \cite[Theorem~4]{GITEigen2}.
\end{proof}

\paragraphe
Consider now the projective space $D=\PP(\Rep(Q,\alpha))$ endowed with the 
$\GL(\alpha)$-action.
We are now interested in the GIT-cone $\ac^{\GL(\alpha)}(D)$.

We have the following obvious

\begin{lemma}
\label{lem:PicGD}
  The restriction map $\rho_D\,:\,\Pic^{\GL(\alpha)}(X)\longto\Pic^{\GL(\alpha)}(D)$ is an isomorphism of groups.
Moreover, $\rho_D(\Li)$ is ample if and only if $\Li$ is.
\end{lemma}

Lemma~\ref{lem:PicGD} allows to embed $\ac^{\GL(\alpha)}(D)$ in $\QQ\times\Gamma_\QQ$. 
Set $\Pol(D,Q,\alpha)=\ac^{\GL(\alpha)}(D)\cap \{1\}\times\Gamma_\QQ$.
Obviously,   $\Pol(D,Q,\alpha)$ is a rational polytope in $\Gamma_\QQ$.
\\

Via the identification of Lemma~\ref{lem:PicGD}, the relation between  
$\Pol(D,Q,\alpha)$ and $\Pol(Q,\alpha)$ is as follows:

\begin{prop}
\label{prop:PDalpha}
The polytope $\Pol(Q,\alpha)$ is the convex hull of $0$ and
$\Pol(D,Q,\alpha)$.

We assume moreover that $Q$ has no cycle (even non oriented). Then, 
$\Pol(D,Q,\alpha)$ is a face of $\Pol(Q,\alpha)$.
In particular, $\Pol(D,Q,\alpha)$ is an affine section of $\Sigma(Q,\alpha)$.
\end{prop}

\begin{proof}
Let $\sigma\in\Gamma_\QQ$.
It is clear that $\sigma\in\Pol(D,Q,\alpha)$  if and only if 
$X^{\rm ss}(\sigma)$ intersects $D$,  if and only if
$X^{\rm ss}(\sigma)$ is not contained in $X^{\rm ss}(0)$.
By \cite{GeomDedic}, this is equivalent to the fact that 
the closure of the GIT-class of $\sigma$ does not contain $0$.
In particular, all the vertices of $\Pol(Q,\alpha)$ excepted $0$ belong 
to $\Pol(D,Q,\alpha)$; the first assertion follows.

With the additional assumption, one can easily construct a central one
parameter subgroup $\lambda(t)$ of $G(\alpha)$ which acts on 
each $\Hom(V(ia),V(ta))$ (for $a\in Q_1$) by multiplication by $t$. 
Then, $(D,\lambda)$ is a well covering pair; this implies that
$\Pol(D,Q,\alpha)$ is a face of $\Pol(Q,\alpha)$. 
\end{proof}\\

\paragraphe
Propositions~\ref{prop:0sommet} and \ref{prop:PDalpha} proves that when $Q$ has no cycle,
the descriptions of the three cones are equivalent. From now on, we are mainly 
interested in $\Sigma(Q,\alpha)$ viewed as the cone generated by $\Pol(Q,\alpha)$; 
that is to the faces of $\Pol(Q,\alpha)$ containing the vertex $0$:

\begin{lemma}
\label{lem:C0}
Let $(C,\lambda)$ be a dominant pair. 
Then, $\Face(C)$ contains $0$ if and only if $C$ contains $0$.  
\end{lemma}

\begin{proof}
  The point is that $\{0\}$ is the only closed orbit in $X^{\rm ss}(0)$.
Actually, if $\Face(C)$ contains $0$, $C$ has to contains $0$ by 
\cite[Theorem~5]{GITEigen}. Conversely, if $C$ contains $0$, $\mu^0(C,\lambda)=0$.
\end{proof}

\subsection{Dominant pairs}

\paragraphe\label{par:lambdaalpha}
Let $\sigma\in\Gamma$ and $\alpha$ a vector dimension. We set:
$$
\sigma(\alpha):=\sum_{s\in Q_0}\sigma(s)\alpha(s).
$$
We consider the one parameter subgroup $\lambda_\alpha$ of $\GL(\alpha)$ acting on $V(s)$ by $t.\Id$ for any 
$s\in Q_0$: $\sigma(\alpha)$ is simply the composition $\sigma\circ\lambda_\alpha$. 
Note that $\lambda_\alpha$ acts trivially on $\Rep(Q,\alpha)$. 
This implies that $\Pol(Q,\alpha)$ is contained in the hyperplane $\Hcal(\alpha)$ consisting of the $\sigma$'s
such that $\sigma(\alpha)=0$.\\
 
\begin{defin}
  The dimension vector $\alpha$ is called a {\it rational Schur root} if $\Pol(Q,\alpha)$ or equivalently $\Sigma(Q,\alpha)$
has non empty interior in $\Hcal(\alpha)$.

If there exists  $R\in\Rep(Q,\alpha)$ whose the stabilizer in $\GL(\alpha)$ has dimension one, $\alpha$ is said to be
a {\it Schur root}.
\end{defin}

The second notion is very classical (see \cite{Kac:carquois2}) and the first one very natural in our context. 
We will explain the relation between these two notions in Paragraph~\ref{par:ratSchur}.\\

\paragraphe {\bf Decompositions of dimension vectors.}
Let $\alpha$ be a vector dimension of $Q$.\\

\begin{defin}
A {\it $\ZZ$-decomposition} of $\alpha$ is a family of dimension vectors $\alpha_i$ indexed by $\ZZ$
such that $\alpha_i=0$ with finitely many exceptions and $\alpha=\sum_{i\in\ZZ}\alpha_i$.
An {\it ordered decomposition} of $\alpha$, is a sequence $(\beta_1,\cdots,\beta_s)$ of
non-zero vector dimensions such that $\alpha=\beta_1+\cdots+\beta_s$.
We denote the decomposition by $\alpha=\beta_1\pplus\cdots\pplus\beta_s$.
\end{defin}

\paragraphe \label{par:defeta}
Let $\lambda$ be a one parameter subgroups of $\GL(\alpha)$.
For any $i\in \ZZ$ and $s\in Q_0$, we set $V_i(s)=\{v\in V(s)\,|\,\lambda(t)v=t^iv\}$ and 
$\alpha_i(s)=\dim V_i(s)$.
Obviously, $\alpha=\sum_{i\in\ZZ}\alpha_i$ form a $\ZZ$-decomposition of $\alpha$ which 
determines $\lambda$ up to conjugacy.

The parabolic subgroup $P(\lambda)$ of $\GL(\alpha)$ associated to $\lambda$ is the set of 
$(g(s))_{s\in Q_0}$ such that for all $i\in\ZZ$ and $s\in Q_0$ we have 
$g(s)(V_i(s))\subset \oplus_{j\leq i}V_j(s)$.

Now, $\Rep(Q,\alpha)^\lambda$ is the set of the $(u(a))_{a\in Q_1}$'s such that 
for any $a\in Q_1$ and for any $i\in\ZZ$, $u(a)(V_i(ia))\subset V_i(ta)$. 
It is  isomorphic to $\prod_i\Rep(Q,\alpha_i)$.
In particular, the irreducible component $C$ of $X^\lambda$ containing $0$ is isomorphic
to $\PP(\prod_i\Rep(Q,\alpha_i)\oplus\CC)$.

Moreover, $C^+\cap \Rep(Q,\alpha)$ is the set of the $(u(a))_{a\in Q_1}$'s such that 
for any $a\in Q_1$ and for any $i\in\ZZ$, $u(a)(V_i(ia))\subset \oplus_{j\leq i}V_j(ta)$.

Consider the morphism $\eta_\lambda\,:\,G\times_{P(\lambda)}C^+\longto \Rep(Q,\alpha)$.
Note that, $P(\lambda)$, $C$ and $C^+$ only depend (up to conjugacy) on the ordered decomposition of 
$\alpha$ induced by the $\ZZ$-decomposition $\sum_i\alpha_i$ in an obvious way. From now on, we will 
consider the
{\it map $\eta_{\beta_1\pplus\cdots\pplus\beta_s}$} associated to the ordered decomposition of $\alpha$; 
it is defined up to conjugacy.
We will say that the ordered decomposition is dominant respectively birational if 
$\eta_{\beta_1\pplus\cdots\pplus\beta_s}$ is. We will say that the decomposition is {\it well covering}
if $(C,\lambda)$ is.
If the decomposition is dominant, we will denote by 
$\Face(\beta_1\pplus\cdots\pplus\beta_s)$ the corresponding face of $\Pol(Q,\alpha)$.

\begin{lemma}
\label{lem:descFdecomp}
Let $\beta=\beta_1\pplus\cdots\pplus\beta_s$ be a dominant ordered decomposition.
Then, $\Face(\beta_1\pplus\cdots\pplus\beta_s)$ is generated by the face 
$\Hcal(\beta_1)\cap\cdots\cap\Hcal(\beta_s)\cap\Sigma(Q,\beta)$ of $\Sigma(Q,\beta)$.
\end{lemma}

\begin{proof}
  Let $(C,\lambda)$ be a dominant pair associated to  $\beta=\beta_1\pplus\cdots\pplus\beta_s$.
Let us fix $\underline V=(V(s))_{s\in Q_0}$ of dimension $\beta$.
Let $\underline V=\underline V_1\oplus\underline V_s$ be a decomposition such that $\underline V_i$ 
has dimension $\beta_i$.
The torus $(\CC^*)^s$ acts on $\underline V$ as follows; the $i^{\rm th}$ component acts by homothetie
on $\underline V_i$.
The induced action of $(\CC^*)^s$ on $C$ is trivial. 
So, $(\CC^*)^s$ has to acts trivially on any point in $\Face(C)$; it follows that,
$\Face(C)$ is contained in $\Hcal(\beta_1)\cap\cdots\cap\Hcal(\beta_s)$.

Conversely, let $\sigma\in \Hcal(\beta_1)\cap\cdots\cap\Hcal(\beta_s)\cap\Sigma(Q,\beta)$.
Since $(C,\lambda)$ is covering, $X^{\rm ss}(\Li(1,\sigma))$ intersects $C^+$.
But, since $\sigma\in \Hcal(\beta_1)\cap\cdots\cap\Hcal(\beta_s)$, $\lambda$ acts trivially 
on $\Li(1,\sigma)_{|C}$. By \cite[Lemma~4]{GITEigen}, this implies that $X^{\rm ss}(\Li(1,\sigma))$ 
intersects $C$.
\end{proof}\\

\paragraphe
Let $\alpha,\beta\in\NN^{Q_0}$.
Following Derksen-Schofield-Weyman (see~\cite{DSW:nbsubrep}), we define $\alpha\circ\beta$ to be the 
number of $\alpha$-dimensional subrepresentation of a general representation of dimension $\alpha+\beta$ 
if it is finite, and $0$ otherwise.

We now recall a description of well covering ordered decomposition from \cite{multi}:
\begin{prop}
\label{prop:wcod}
The ordered decomposition $\beta= \beta_1\pplus\cdots\pplus\beta_s$ is well covering if and only if
$$\forall i<j\ \ \beta_i\circ \beta_j=1.
$$
\end{prop}

\paragraphe
\label{par:ratSchur}
We can now explain the name ``rational Schur root''.
Let us first reprove two well known lemmas:

\begin{lemma}
\label{lem:alphabetagroup}
If $\alpha\circ\beta\neq 0$ and $\alpha\circ\gamma\neq 0$ then 
$\alpha\circ(\beta+\gamma)\neq 0$.  
\end{lemma}

\begin{proof}
In \cite{DSW:nbsubrep}, Derksen-Schofield-Weyman proved that $\alpha\circ\beta$ is the dimension of 
$\CC[\Rep(Q,\alpha)]_\sigma$ for well chosen weight $\sigma$. With this characterization, the lemma 
just follows from the fact that $\CC[\Rep(Q,\alpha)]^{\SL(\alpha)}$ is an algebra.
In this work, $\alpha\circ\beta$ is always the degree a map $\eta$; and, we include a proof using this
point of view.

Consider a pair $(C,\lambda)$ (resp. $(C',\lambda'))$ associated to the ordered decomposition
$\alpha\pplus\beta$ (resp. $\alpha\circ\gamma$) in $\Rep(Q,\alpha+\beta)$ and $\Rep(Q,\alpha+\gamma)$.
Since $\alpha\circ\beta\neq 0$, $\eta_{\alpha\pplus\beta}$ is generically finite.
Moreover, by \cite[Lemma~]{multi}, $\lambda$ acts trivially on the restriction to $C$ of the determinant
bundle of $\eta$. It follows that for general $x\in C=\Rep(Q,\alpha)\oplus\Rep(Q,\beta)$, the differential of 
$\eta_{\alpha\pplus\beta}$ at $x$ is an isomorphism.
In the same way, the differential of $\eta_{\alpha\pplus\gamma}$ is an isomorphism for $x'$ general in
$\Rep(Q,\alpha)\oplus\Rep(Q,\beta)$. 
A direct computation implies now that $\eta_{\alpha\pplus(\beta+\gamma)}$ is an isomorphism 
for $y$ general in  $\Rep(Q,\alpha)\oplus\Rep(Q,\beta)\oplus\Rep(Q,\gamma)\subset\Rep(\alpha+\beta+\gamma)$.
In particular, $\alpha\circ(\beta+\gamma)\neq 0$. 
\end{proof}

\begin{lemma}
  \label{lem:Sigmafirstdesc}
We have:
$$
\Sigma(Q,\beta)=\{\sigma\in\Gamma\,:\,\sigma(\beta)=0{\rm\ and\ }\sigma(\alpha)\leq 0\ \forall\alpha
{\rm\ s.t.\ }\alpha\circ(\beta-\alpha)\neq 0\}.
$$
\end{lemma}

\begin{proof}
Let $\sigma\in\Sigma(Q,\beta)$.
We already saw that $\sigma(\beta)=0$. 
Let $\alpha$ be such that $\alpha\circ(\beta-\alpha)\neq 0$.
Since $\eta_{\alpha\pplus(\beta-\alpha)}$ is dominant, $\sigma(\alpha)\leq 0$.

The converse inclusion is a direct consequence of \cite{King:moduli} 
(see~\cite[Remark~5]{DW:saturation}).
\end{proof}

Here, comes an easy variant of Derksen-Schofield-Weyman's saturation Theorem:

\begin{lemma}
  \label{lem:Sigmakbeta}
We have:
$$
\Sigma(Q,k\beta)=\Sigma(Q,\beta).
$$
\end{lemma}

\begin{proof}
The inclusion   $\Sigma(Q,k\beta)\subset \Sigma(Q,\beta)$ is a direct consequence of 
Lemmas~\ref{lem:alphabetagroup} and \ref{lem:Sigmafirstdesc}.

The converse inclusion follows from Derksen-Weyman's Reciprocity Property 
(see \cite[Corollary~1]{DW:saturation}). We include here a simpler proof.
Let $(V(s))_{s\in Q_0}$ be vector spaces of dimension vector $\beta$.  
Consider $(\Hom(\CC^k,V(s))_{s\in Q_0}$ as a $Q_0$ family of vector spaces of dimension vector $k\beta$.
Then, for the natural inclusion  $\Rep(Q,\beta)\subset \Rep(Q,k\beta)$,  $\Rep(Q,\beta)$ is the fix point set of 
$H=(\GL_k)^{Q_0}\subset\GL(k\beta)$. Moreover, the centralizer of $(\GL_k)^{Q_0}$ in $\GL(k\beta)$ is isomorphic
to $\GL(\beta)$. By a Luna's Theorem (see \cite{Luna:adh}), for any linearized ample line bundle a point 
$x\in \Rep(Q,\beta)$ is semistable for $\Li$ and the action of $\GL(\beta)$ if and only if it is for the 
action of $\GL(k\beta)$.
It follows that $\Pol(Q,\beta)\subset\Pol(Q,k\beta)$.
The lemma is proved.
\end{proof}

\begin{prop}
  \label{prop:}
A vector dimension $\alpha$ is a rational Schur root if and only if it is positively proportional 
to a Schur root.
\end{prop}

\begin{proof}
The Ringle form is denoted by  $\langle\cdot,\cdot\rangle$.
Let $\beta$ be a Schur root.
By \cite[Theorem~6.1]{Scho:genrep}, $X$ contains stable points for the action of 
$\GL(\beta)/{\rm Im}(\lambda_0)$
and the line bundle $\Li(1,\langle\beta,\cdot\rangle-\langle\cdot,\beta\rangle)$.
It follows that $\Sigma(Q,\beta)$ has non empty interior in $\Hcal(\beta)$.
By Lemma~\ref{lem:Sigmakbeta}, $k\beta$ is a rational Schur root for any positive integer $k$.

Conversely, let $\beta$ be a rational Schur root. Let $d$ denote the gcd of the $\beta(s)$ for $s\in Q_0$.
By Lemma~\ref{lem:Sigmakbeta}, $\overline{\beta}=\beta/d$ is a rational Schur root.
Consider the canonical decomposition $\overline{\beta}=\beta_1+\cdots+\beta_s$ of $\beta$ 
(see \cite{Kac:carquois2}). 
Then, $\Sigma(Q,\overline{\beta})$ is contained in $\Hcal(\beta_1)\cap\cdots\cap\Hcal(\beta_s)$.
Since  $\Sigma(Q,\overline{\beta})$ spans the hyperplane $\Hcal(\beta)$, it follows that 
$\Hcal(\beta)=\Hcal(\beta_1)=\cdots=\Hcal(\beta_s)$.
So, the $\beta_i$'s are proportional; since, $\overline{\beta}$ is indivisible, it follows that $s=1$ and
that $\overline{\beta}$ is a Schur root.
\end{proof}

\subsection{Derksen-Weyman's Theorem}

\paragraphe
The ordered decomposition $\beta= \beta_1\pplus\cdots\pplus\beta_s$ is called an 
{\it ordered decomposition by rational Schur roots} if $\beta_1,\cdots,\beta_s$ are rational 
Schur roots.
To any such decomposition we associate the (unordered) set $\{\beta_1,\cdots,\beta_s\}\subset \NN^{Q_0}$.
Let $\Wcal_s(\beta)$ denote the set of subsets  obtained in such a way from
well covering ordered decomposition by $s$ rational Schur roots.

We can now state and prove Derksen-Weyman's Theorem:

\begin{theo}
\label{th:main}
Let $\beta$ be a vector dimension. Let $d$ denote the dimension of $\Sigma(Q,\beta)$ and $n$ the cardinality of
$Q_0$. 
For any $s=n-d,\cdots,0$, the map 
$$
\begin{array}{cccc}
  \Theta\ :&\Wcal_s(\beta)&\longto&\{{\rm faces\ of\ }\Sigma(Q,\beta){\rm\ of\ codimension\ }d\}\\
&\{\beta_1,\cdots,\beta_s\}&\longmapsto&\Hcal(\beta_1)\cap\cdots\cap\Hcal(\beta_s)\cap\Sigma(Q,\beta),
\end{array}
$$  
is a bijection.
Moreover, the family $(\beta_1,\cdots,\beta_s)$ is linearly independent.
\end{theo}

\begin{proof}
In this proof, we prefer to consider the faces of $\Pol(Q,\alpha)$ containing $0$ rather than faces of 
$\Sigma(Q,\alpha)$. By Proposition~\ref{prop:0sommet}, this is equivalent.

Let $\beta= \beta_1\pplus\cdots\pplus\beta_s$ be a well covering ordered decomposition by rational Schur roots.
Then,  by Lemma~\ref{lem:descFdecomp}, $\Hcal(\beta_1)\cap\cdots\cap\Hcal(\beta_s)\cap\Sigma(Q,\beta)$ is the face
of $\Sigma(Q,\beta)$ corresponding to $\Face(\beta_1\pplus\cdots\pplus\beta_s)$. 
Since the $\beta_i$'s are rational Schur roots, Theorem~\ref{th:wc2face} shows that 
$\Face(\beta_1\pplus\cdots\pplus\beta_s)$ has codimension $s$.
Let us recall that $\ac^{\GL(\beta)}(X)=\tc^{\GL(\beta)}(X)$.
This proves that $\Theta$ is well defined.\\

Let us fix a face $\Face$ of $\Pol(Q,\alpha)$ of codimension $d$.
By Theorem~\ref{th:faceL}, the exists an open subset $U$ in $\Pic^G(X)^+_\QQ-\ac^G(X)$ such that $\Face=\Face(\Li)$ 
for all $\Li\in U$.
Let $(C,\lambda)$ be a well covering pair associated to a line bundle $\Li\in U$: by Lemma~\ref{lem:C0}, $C$ contains $0$. 
Let $\beta=\beta_1\pplus\cdots\pplus\beta_s$ be the ordered decomposition associated to $\lambda$.
By Paragraph~\ref{par:defeta}, $\eta_{(C,\lambda)}=\eta_{\beta_1\pplus\cdots\pplus\beta_s}$.
  
We claim that the $\beta_i$'s are rational Schur roots.
Let us fix $i\in\{1,\cdots,s\}$. Let $\lambda_{\beta_i}$ be the central one parameter subgroup 
of $\GL(\beta_i)$ defined in Paragraph~\ref{par:lambdaalpha}; and, $S_i$ be the codimension one subtorus of 
the center of  $\GL(\beta_i)$ such that $Y(S_i)$ is orthogonal to $\lambda_{\beta_i}$.
Consider the subgroup $H_i$ of $\GL(\beta_i)$ generated by the $S_i$ and $\SL(\beta_i)$.
We embed $\PP(\Rep(Q,\beta_i)\oplus \CC)$ in $X$ in an obvious way and consider the restriction morphism:
$$
p_i\,:\,\Pic^G(X)_\QQ\longto\Pic^{H_i}(\PP(\Rep(Q,\beta_i)\oplus \CC))_\QQ.
$$ 
By construction, the restriction of $p_i$ to $\Hcal(\beta_i)$ is surjective.
Moreover, by \cite[Lemma~11]{GITEigen}, $p_i(U)$ is contained in $\Pol(Q,\beta_i,H_i)$.
Since $p_i$ is an open map, this implies that $\Pol(Q,\beta_i)$ has codimension one in $X(\GL(\beta_i))_\QQ$.
So, the $\beta_i$'s are rational Schur roots and $\Theta$ is surjective. \\ 

Let $\beta= \beta_1\pplus\cdots\pplus\beta_s$ be any well covering ordered decomposition by rational 
Schur roots.
By Theorem~\ref{th:wc2face}, the intersection $\Hcal(\beta_1)\cap\cdots\cap\Hcal(\beta_s)$ has codimension $s$.
This means that the $\beta_i$ are linearly independent.

Let us fix $\sigma\in\Gamma_\QQ$ such that $\Li(1,\sigma)$ belongs to the relative interior of
$\Face:=\Face( \beta_1\pplus\cdots\pplus\beta_s)$.
Since, $\beta_i$ are rational Schur roots, Theorem~\ref{th:wc2face} shows that the codimension of $\Face$ 
equals $s$.
Note that, $\Delta(\Face)$ contains $\Li(1,\sigma)$ in its closure.
We claim that $p_i(\sigma)$ belongs to the relative interior of $\Sigma(Q,\beta_i)$.
Actually, if not, one can find $\sigma_\epsilon$ in $\Delta(\Face)$ such that 
$p_i(\sigma_\epsilon)$ does not belongs to  $\Sigma(Q,\beta_i)$.
Then, the ordered decomposition associated to $\sigma_\epsilon$ contains strictly more than $s$ vector
dimensions. By Theorem~\ref{th:wc2face} this implies that the codimension of $\Face$ is strictly greater 
than $s$; which is a contradiction.\\

We now want to prove the injectivity of $\Theta$. 
Let $\beta= \beta_1\pplus\cdots\pplus\beta_s$ be a well covering ordered decomposition by rational 
Schur roots and $\Face$ be the associated face.
We want to obtain the decomposition of $\beta$ from $\Face$. 
By Proposition~\ref{prop:} and \cite[Theorem~3.2]{Scho:genrep}, the canonical decomposition of 
$\beta_i=a_i\overline{\beta}_i$ for some positive integer $a_i$ and some Schur root $\overline{\beta_i}$.
Set $C=\PP(\oplus_i\Rep(Q,\beta_i)\oplus\CC)$ and 
$C_0=\PP(\oplus_i\Rep(Q,\overline{\beta}_i)^{\oplus a_i}\oplus\CC)$; and fix embeddings $C_0\subset C\subset X$.

Let $\Li:=\Li(1,\sigma)$ be a point in the relative interior of $\Face$.
Let $x$ be a general point in $C_0$. 
Since $p_i(\sigma)$ belongs to the relative interior of $\Sigma(Q,\overline{\beta}_i)$, 
\cite[Theorem~6.1]{Scho:genrep} implies that the orbit of $x$ by the group 
$\prod_i\GL(\overline{\beta}_i)^{a_i}$ is closed in $X^{\rm ss}(\Li)$.
By \cite{Luna:adh}, this implies that $\GL(\beta).x$ is closed in  $X^{\rm ss}(\Li)$.
Conversely, by \cite[Theorem~5]{GITEigen}, any general closed closed orbit in $X^{\rm ss}(\Li)$ intersects 
$C$ and so $C_0$.
This proves that a general closed isotropy in $X^{\rm ss}(\Li)$ contains a general point of $C_0$.
In particular, any point in a general closed isotropy of $X^{\rm ss}(\Li)$ decompose as 
a sum of $a_1$ indecomposable representations of dimension $\overline{\beta}_1$\dots and $a_s$ 
 indecomposable representations of dimension $\overline{\beta}_s$.
Moreover, such a decomposition is unique and the $\overline{\beta}_i$'s are pairwise distinct 
(the family is free). 
The injectivity follows.
\end{proof}

\bibliographystyle{amsalpha}
\bibliography{biblio}

\providecommand{\bysame}{\leavevmode\hbox to3em{\hrulefill}\thinspace}
\providecommand{\MR}{\relax\ifhmode\unskip\space\fi MR }
\providecommand{\MRhref}[2]{%
  \href{http://www.ams.org/mathscinet-getitem?mr=#1}{#2}
}
\providecommand{\href}[2]{#2}
\begin{thebibliography}{DSW07}

\bibitem[BK06]{BK}
Prakash Belkale and Shrawan Kumar, \emph{Eigenvalue problem and a new product
  in cohomology of flag varieties}, Invent. Math. \textbf{166} (2006), no.~1,
  185--228.

\bibitem[BS00]{sjamaar}
Arkady Berenstein and Reyer Sjamaar, \emph{Coadjoint orbits, moment polytopes,
  and the {H}ilbert-{M}umford criterion}, J. Amer. Math. Soc. \textbf{13}
  (2000), no.~2, 433--466 (electronic).

\bibitem[DSW07]{DSW:nbsubrep}
Harm Derksen, Aidan Schofield, and Jerzy Weyman, \emph{On the number of
  subrepresentations of a general quiver representation}, J. Lond. Math. Soc.
  (2) \textbf{76} (2007), no.~1, 135--147.

\bibitem[DW00]{DW:saturation}
Harm Derksen and Jerzy Weyman, \emph{Semi-invariants of quivers and saturation
  for {L}ittlewood-{R}ichardson coefficients}, J. Amer. Math. Soc. \textbf{13}
  (2000), no.~3, 467--479 (electronic).

\bibitem[DW06]{DW:comb}
Harm Derksen and Jerzy Weyman, \emph{The combinatorics of quiver
  representations}, 2006, {\tt arXiv.org:math/0608288}.

\bibitem[Kac82]{Kac:carquois2}
V.~G. Kac, \emph{Infinite root systems, representations of graphs and invariant
  theory. {II}}, J. Algebra \textbf{78} (1982), no.~1, 141--162.

\bibitem[Kin94]{King:moduli}
A.~D. King, \emph{Moduli of representations of finite-dimensional algebras},
  Quart. J. Math. Oxford Ser. (2) \textbf{45} (1994), no.~180, 515--530.

\bibitem[Lun75]{Luna:adh}
D.~Luna, \emph{Adh\'erences d'orbite et invariants}, Invent. Math. \textbf{29}
  (1975), no.~3, 231--238.

\bibitem[Res00]{GeomDedic}
N.~Ressayre, \emph{The {GIT}-equivalence for {$G$}-line bundles}, Geom.
  Dedicata \textbf{81} (2000), no.~1-3, 295--324.

\bibitem[Res07]{GITEigen}
Nicolas Ressayre, \emph{Geometric invariant theory and generalized eigenvalue
  problem}, Preprint (2007), no.~arXiv:0704.2127, 1--45.

\bibitem[Res08a]{GITEigen2}
\bysame, \emph{Geometric invariant theory and generalized eigenvalue problem
  ii}, Preprint (2008), 1--25.

\bibitem[Res08b]{multi}
\bysame, \emph{Multiplicative formulas in cohomology of g/p and in quiver
  representations}, Preprint (2008), no.~arXiv:0812.2122, 1--20.

\bibitem[Sch92]{Scho:genrep}
Aidan Schofield, \emph{General representations of quivers}, Proc. London Math.
  Soc. (3) \textbf{65} (1992), no.~1, 46--64.

\end{thebibliography}

\begin{center}
  -\hspace{1em}$\diamondsuit$\hspace{1em}-
\end{center}

\end{document}